\documentclass{article}                                                  % if you need a4paper
\usepackage{algorithm, algorithmic, setspace,amsthm}
\usepackage{graphicx} % for pdf, bitmapped graphics files
\usepackage{amsmath} % assumes amsmath package installed
\usepackage{amssymb}  % assumes amsmath package installed
\usepackage{amsfonts}
\usepackage{color}
\usepackage{subfig}
\usepackage[english]{babel}
\usepackage[latin1]{inputenc}
\usepackage{siunitx}
\usepackage{steinmetz}
\newtheorem{remark}{Remark}

\newtheorem{theorem}{Theorem}

\newtheorem{proposition}{Proposition}
\newcommand{\sto}{ \rm{s.t.} }

\newcommand{\norm}[1]{\left\lVert #1 \right\rVert}

\newcommand{\Lag}{\mathcal{L}}
\newcommand{\zstar}{z^{\star}}
\newcommand{\xstar}{x^{\star}}
\newcommand{\lstar}{\lambda^{\star}}
\newcommand{\xstart}{x^{\star\top}}
\newcommand{\lstart}{\lambda^{\star\top}}

\newcommand{\R}{\mathbb{R}}
\newcommand{\ki}{K_i}
\newcommand{\kp}{K_p}
\usepackage{soul}

\title{\LARGE \bf A feedback control approach to convex optimization with inequality constraints}

\author{V. Cerone, S. M. Fosson, S. Pirrera, D. Regruto\thanks{Corresponding author. The authors are with the Dipartimento di Automatica e Informatica,
Politecnico di Torino, Corso Duca degli Abruzzi 24, 10129 Torino, Italy;
    e-mail: diego.regruto@polito.it.
}
}

\begin{document}
\maketitle
\thispagestyle{empty}
\pagestyle{empty}

\begin{abstract}
We propose a novel continuous-time algorithm for inequality-constrained convex optimization inspired by proportional-integral control. Unlike the popular primal-dual gradient dynamics, our method includes a proportional term to control the primal variable through the Lagrange multipliers. This approach has both theoretical and practical advantages. On the one hand, it simplifies the proof of the exponential convergence in the case of smooth, strongly convex problems, with a more straightforward assessment of the convergence rate concerning prior literature. On the other hand, through several examples, we show that the proposed algorithm converges faster than primal-dual gradient dynamics. This paper aims to illustrate these points by thoroughly analyzing the algorithm convergence and discussing some numerical simulations.
\end{abstract}

\section{Introduction}\label{sec:IN}
Primal-dual gradient dynamics (PDGD) is a well-established continuous-time algorithm that solves constrained optimization problems. Introduced in \cite{kos56,arr58}, it consists of a primal-descent, dual-ascent gradient method achieving the saddle point of the Lagrangian of the problem.

In the last years, we have witnessed a renewed interest in PDGD thanks to its effectiveness in several engineering applications and control problems, e.g., game theory \cite{gha13}, power systems \cite{zha14,na16} and model predictive control \cite{ski23}. A particular focus is on its use in distributed optimization; see, e.g., \cite{fei10, gha13, cor19, shi23}. Gradient-based algorithms are notably well-suited for implementation over decentralized networks.

In the recent literature, several works have addressed the study of the stability and convergence of PDGD. This algorithm is globally exponentially convergent for smooth, strongly convex problems (see, e.g., \cite{qu19,alg20,che16}) and for problems that combine strongly and non-strongly convex terms in \cite{du19,che20}. In \cite{che17} and \cite{che18}, the authors study the asymptotic convergence for general saddle functions not directly related to constrained optimization. In \cite{din19} and \cite{dhi19}, the analysis envisages also nonsmooth composite optimization problems.
Among the mentioned works, \cite{alg20} and \cite{che20} consider equality-constrained problems, while \cite{che16} and \cite{qu19} also consider inequality constraints.

This paper proposes a novel continuous-time algorithm for smooth, strongly convex problems with inequality constraints. By starting from the definition of a suitable augmented Lagrangian, the key idea is to control the dynamics of the primal variable through the Lagrange multipliers of the problem by implementing a feedback control method inspired by proportional-integral (PI) control. The contribution of the paper is twofold. On the one hand, we prove the exponential convergence of the proposed method for strongly convex functions. On the other hand, we show its practical effectiveness through numerical simulations. In particular, we analyze its behavior when compared to PDGD.

This work partially extends the framework proposed in \cite{cmo24}, where we develop a feedback control approach for equality-constrained problems, specializing in PI control and feedback linearization. In this paper, we retrieve the key ideas of the PI control algorithm proposed in \cite{cmo24}, and we develop a novel PI approach for the case of inequality constraints. In particular, this extension requires a novel convergence analysis starting from a peculiar augmented Lagrangian.

We organize the paper as follows. Section \ref{sec:PS} states the problem and reviews its solution through PDGD. In Section \ref{sec:PA}, we develop the proposed algorithm while we study its convergence in \ref{sec:CA}. Section \ref{sec:NR} shows the effectiveness of the proposed method through numerical experiments, with particular attention to the convergence speed. Finally, Section \ref{sec:CON} concludes the paper.
\section{Problem statement and related work}\label{sec:PS}
Let $f: \R^n \to \R$ be a smooth, strongly convex function. We consider the constrained optimization problem
\begin{equation} \begin{aligned}
    \min_{x \in \R^n}& \, f(x) \\ &\sto \\ &h(x)=Cx-d\leq 0
    \label{opt}
\end{aligned} \end{equation}
where $C\in\R^{m,n}$, $d\in\R^m$ and ``$\leq$'' denotes the componentwise inequality. %Let $d_j$ and $C_j$ be the $j$-th component of $d$ and the $j$-th row of $C$, respectively.
%We call $\Omega \doteq \{ x | h(x) = 0 \}$ the feasible set; any $x \in \Omega$ is a feasible point.
%
{ As in \cite{qu19},} we consider the following augmented Lagrangian
\begin{equation}\label{lag}
\Lag(x,\lambda)=f(x)+g(x,\lambda)
\end{equation}
where we define { $g:\R^{n} \times \R^m \to \R$ }as
\begin{equation}\label{g}
 \begin{split}
  &g(x,\lambda)=\sum_{j=1}^m  g_j(x,\lambda_j)~~ \text{ with }\\
   &g_j(x,\lambda_j)=\\&=\left\{\begin{aligned}
                             &\lambda_j h_j(x)+\frac{\rho}{2}h_j^2(x) ~~~~~\text{ if } h_j(x)\geq -\frac{\lambda_j}{\rho};\\
                             &-\frac{1}{2\rho}\lambda_j^2~~~~~~~~~~~~~~~~~\text{ otherwise}
                            \end{aligned}\right.
 \end{split}
\end{equation}
where $\rho>0$ is a design hyperparameter. %{We remark that the selection of $g(x,\lambda)$ in \eqref{g} is standard in the context of inequality-constrained optimization. See, e.g., \cite{ber82,qu19} for further details.}
%% in realtà non è quella standard, vedi  Remark 2 di Qu 2019: questo è anche detto qualche riga sotto

The function $g_j(x,\lambda_j):\R^{n+1}\to\R$ penalizes the constraint violation. We notice that it is continuous and has a continuous gradient for each $(x,\lambda_j)\in\R^{n+1}$. 
As shown in \cite[Eq. (9)]{qu19}, PDGD for problem \eqref{lag}-\eqref{g} corresponds to the dynamic system
\begin{equation}\label{PDGD}
   \begin{split}
    &\dot{x}=-\nabla_x \Lag(x,\lambda)=-\nabla f(x)-\nabla_x g(x,\lambda)\\
    &\dot{\lambda}=\eta \nabla_{\lambda}g(x,\lambda)
   \end{split}
  \end{equation}
  for some $\eta>0$.
Let $J_h(x)\in\R^{m,n}$ be the Jacobian matrix of $h$, i.e., $J_h(x)=C$ in the linear case. Then,
\begin{equation}\label{gradienti}
 \begin{split}
  &\nabla_x g_j(x,\lambda_j)=\max\{\rho h_j(x)+\lambda_j,0\}(J_h(x))_j^\top\\
  &\frac{\partial g_j}{\partial \lambda_j}(x,\lambda_j)=\frac{\max\{\rho h_j(x)+\lambda_j,0\}-\lambda_j}{\rho}.
 \end{split}
\end{equation}
PDGD in \eqref{PDGD} is a system switching between two modes for each $j$. In the first mode, i.e., when
$h_j(x)\geq -\frac{\lambda_j}{\rho}$, the Lagrange multipliers $\lambda(t)$ {\it{control}} the dynamics of the state $x(t)$, in order to push the system towards the feasible set $h(x)\leq 0$. In the second mode, i.e., when the constraint is satisfied, $x(t)$ is not controlled and evolves based on the gradient of $f$; contextually, $\lambda(t)$ converges to zero.

We remark that \eqref{lag} is not the standard Lagrangian for inequality-constrained problems used, e.g., in \cite{fei10,che16}. As noticed in \cite{qu19}, the standard Lagrangian gives rise to a PDGD with a discontinuous projection step, which creates numerical issues when implementing the algorithm. Moreover, using the standard Lagrangian makes the convergence analysis more challenging, and only stability is proven; the authors of \cite{qu19} conjecture that this discontinuous PDGD is not globally exponentially stable.

As noticed in \cite{qu19}, the saddle point of \eqref{lag}-\eqref{g} corresponds to the saddle point of the standard Lagrangian; see, e.g., \cite[Chapter 3]{ber82} for details.

The first-order Karush-Kuhn-Tucker (KKT) conditions for problem \eqref{opt} are as follows, see, e.g., \cite[Section 11.8]{lue_book}:
\begin{theorem}[KKT first-order  conditions]
\label{lagr_th}
Let $\xstar \in \R^n$ be a minimum of $f$ subject to the constraints $h(x)\leq 0$. Then, there exists a unique  $\lstar \in \R^m$ such that
\begin{equation} \label{FO}
\begin{split}
   &\nabla f(\xstar) + J_h(\xstar)^{\top}\lstar    = 0\\
   &\lstart h(\xstar)=0\\
   &\lstar\geq 0.
   \end{split}
\end{equation}
%where $J_h(x)\in\R^{m,n}$ is the Jacobian matrix of $h$ evaluated in $x$.
% Moreover, if $f,h$ are twice differentiable,
% \begin{equation}
%     y^\top \left(\hess{f}{x}(\xstar) + \sum_{i=1}^{m} \lambda_i^* \hess{h_i}{x}(\xstar)\right) y \geq 0
%     \label{lagr_necess2}
% \end{equation}
% for all $y \in \left\{ y | \left(\nabla{h_i}{x}(\xstar)\right)^\top y = 0,\, i=1,\dots,m \right\}$.
\end{theorem}
Proposition 1 in \cite{qu19} states that PDGD has a unique equilibrium point $(\xstar,\lstar)$, and it satisfies the first-order KKT conditions for problem \eqref{opt}. Furthermore, Theorem 2 in \cite{qu19} proves that PDGD is globally exponentially convergent. The computed convergence rate $\tau_{ineq}$ depends on several constants arising from repeated upper/lower bounding of the eigenvalues of the matrices involved in the proof; therefore, it is difficult to explicitly assess $\tau_{ineq}$ from the formula given in  Theorem 2 in \cite{qu19}. The proof is rather technical because it needs a { non-diagonal quadratic Lyapunov function}; see \cite[Eq. (10)]{qu19}.

As to strongly convex problems with $h(x)=0$, in \cite{cmo24}, we show that PDGD corresponds to an integral control system that regulates $h(x)$ to zero, based on standard Lagrangian for equality constraints problem, i.e., $f(x)+g_{eq}(x,\lambda)$ with $g_{eq}(x,\lambda)=\lambda^\top h(x)$. In this system, the Lagrange multipliers $\lambda$ play the role of the control input. In \cite{cmo24}, we design a PI control, which usually has a faster convergence rate than PDGD. This PI control system reads as follows:
\begin{equation}\label{PI_EQ}
   \begin{split}
    &\dot{x}=-\nabla_x \Lag(x,\lambda)=-\nabla f(x)-\nabla_x g_{eq}(x,\lambda)\\
     &\dot{\lambda}=\ki \nabla_{\lambda}g_{eq}(x,\lambda)+\kp\frac{d}{dt}\nabla_{\lambda}g_{eq}(x,\lambda)   \end{split}
  \end{equation}
   where $\ki$ and $\kp$ are design hyperparameters.
Since $\nabla_{\lambda}g_{eq}(x,\lambda)=h(x)$, this a PI control regulating $h(x)$ to zero.
\section{Proposed approach}\label{sec:PA}
This section proposes a novel PI control approach to solve \eqref{opt}. A natural extension of \eqref{PI_EQ} to the case of inequality constraints would be the application of \eqref{PI_EQ} to the augmented Lagrangian \eqref{lag} with \eqref{g}, i.e.,
\begin{equation}\label{PI_NO}
   \begin{split}
    &\dot{x}=-\nabla_x \Lag(x,\lambda)=-\nabla f(x)-\nabla_x g(x,\lambda)\\
     &\dot{\lambda}=\ki \nabla_{\lambda}g(x,\lambda)+\kp\frac{d}{dt}\nabla_{\lambda}g(x,\lambda)
   \end{split}
  \end{equation}
{We can interpret the dynamic system \eqref{PI_NO} as a PI controller that pushes $h_j(x)$ towards zero whenever $h_j(x)\geq0$ and provides no control action when $h_j(x)\leq 0$, letting $\lambda_j$ converge to zero.}

%In particular, we notice that $\frac{d}{dt}\nabla_{\lambda}g_j(x,\lambda_j)=\mathbf{1}(h_j(x)\geq -\frac{\lambda_j}{\rho})(J_h(x))_j^\top \dot{x}+\left(\mathbf{1}(h_j(x)\geq -\frac{\lambda_j}{\rho})-1\right)\dot{\lambda}$ where $\mathbf{1}$ is the indicator function.

Even if \eqref{PI_NO} is a natural extension of \eqref{PI_EQ}, proving its {exponential} convergence is challenging. {In particular, we notice that $\frac{d}{dt}\nabla_{\lambda}g_j(x,\lambda_j)$ is discontinuous, which also means the right-hand side of the differential equation describing the closed-loop system is discontinuous and, in turn, potentially, the solution may be non-unique.} For this motivation, we modify \eqref{PI_NO} as follows:
 \begin{equation}\label{PI}
   \begin{split}
    &\dot{x}=-\nabla_x \Lag(x,\lambda)=-\nabla f(x)-\nabla_x g(x,\lambda)\\
    &\dot{\lambda}=\ki \nabla_{\lambda}g(x,\lambda)+\kp J_h(x) \dot{x}
   \end{split}
  \end{equation}
where the involved gradients are explicitly computed in \eqref{gradienti}. We replace $\frac{d}{dt}\nabla_{\lambda}g(x,\lambda)$ by $J_h(x)\dot{x}$, {which we can interpret as taking a continuous approximation of \eqref{PI_NO}.} %We notice that $\frac{d}{dt}\nabla_{\lambda}g_j(x,\lambda_j)=\mathbf{1}(h_j(x)\geq -\frac{\lambda_j}{\rho})(J_h(x))_j^\top \dot{x}+\left(\mathbf{1}(h_j(x)\geq -\frac{\lambda_j}{\rho})-1\right)\dot{\lambda}$ where $\mathbf{1}$ is the indicator function. 
In \eqref{PI_NO}, $\dot{\lambda}$ does not depend on $\dot{x}$ when the constraints are satisfied. In contrast, in \eqref{PI}, the dependence on $\dot{x}$ is always present. %tieThis simplifies the analysis.

As observed for \eqref{PDGD}, PI in \eqref{PI} is a system that switches between two modes for each $j$. In the first mode, i.e., when
$h_j(x)\geq -\frac{\lambda_j}{\rho}$, the Lagrange multipliers $\lambda(t)$ {\it{control}} the dynamics of the state $x(t)$ to achieve $h(x)\leq 0$. {In the second mode, $\lambda(t)$ does not control $x(t)$}. On the other hand, $\lambda(t)$ still depends on $\dot{x}(t)$.

The difference between the proposed approach and PDGD in \eqref{PDGD} and \eqref{PI} is the presence of the additional term $\kp J_h(x) \dot{x}$ in the dynamics of $\lambda$.
To understand the rationale of this term, we go through a feedback control interpretation, as introduced in \cite{cmo24} for the equality-constrained case.
By extending this feedback control framework to the case $h(x)\leq 0$, we can interpret PDGD as an algorithm with integral control on $\lambda_j$ representing a non-satisfied {constraint}. On the other hand, in the presence of a satisfied constraint, we do not control the system through $\lambda$.

In \eqref{PI}, we modify the dynamics of $\lambda$ by adding $\kp J_h(x) \dot{x}$. In the following, we show the benefits of this adjustment in terms of convergence rate.

We remark that \eqref{PI} is not the direct extension of the PI method for $h(x)=0$ reported in \eqref{PI_EQ}. Although feasible, the use of \eqref{PI_EQ} produces {a non-causal system with switched dynamics, creating issues} in the proof of convergence. In other words, through the term $\kp J_h(x) \dot{x}$, we control $\lambda$ via state feedback even when the constraints are satisfied. This modification enhances convergence and simplifies analysis, as shown in the following sections.

%We can interpret the dynamic system  \eqref{PI} can be interpreted as feedback control system, with state $x(t)$ and control input represented by the Lagrange multipliers $\lambda$. We refer to this approach as to controlled Lagrange multipliers (CLM). Specifically, $\lambda$  is computed by feedbacking the state to the controller, whose control law is defined as in the second equation of \eqref{PI}. Moreover, we can consider $\nabla_{\lambda}g(x,\lambda)$ as the output of the system, that we aim to  regulate to zero. More precisely, for each $j=1,\dots,m$, we want to regulate either $h_j(x)$ to zero or $\lambda_j$ to zero; the latter corrsponds to switching off the control input when not necessary to constrain the solution. Of course one might  simply set $\lambda_j=0$ when the control is not necessary, instead of considering a dynamics convergent to zero. However, a smooth behavior for $\lambda$ facilitates the convergence analysis of the overall system.

\section{Convergence analysis}\label{sec:CA}
In this section, we analyse the convergence of the dynamic system \eqref{PI}.
We define
\begin{equation}
 z(t):=(x(t)^\top,\lambda(t)^\top)^\top
\end{equation}
and
\begin{equation}
 \zstar:=(\xstart,\lstart)^\top
\end{equation}
is the equilibrium point of \eqref{PI}, which corresponds to a saddle point of $\Lag(x,\lambda)$.
The following result holds.
\begin{proposition}
The equilibrium point of \eqref{PI} satisfies the KKT conditions \eqref{FO} for problem \eqref{opt}.
\end{proposition}
\begin{proof}
Since at the equilibrium point the time derivatives are null, i.e., $\dot{\xstar}=\dot{\lstar}=0$, we have $\nabla_{\lambda}g(\xstar,\lstar)=0$, i.e., for each $j=1,\dots,m$,
$$\max\{\rho h_j(\xstar)+\lstar_j,0\}=\lstar_j,$$
which implies $\lstar_j\geq 0$, $h_j(\xstar)\leq 0$ and $(h_j(\xstar))\lstar_j=0$.  Finally,
\begin{equation}
\begin{split}
 \nabla_x g(\xstar,\lstar)&=\sum_{j=1}^m \max\{\rho h_j(\xstar)+\lstar_j,0\}C_j^\top\\
 &=\sum_{j=1}^m \lstar_j C_j^{\top}=C^{\top}\lstar=J_h(x)^\top\lstar
 \end{split}
\end{equation}
Therefore,
$\dot{\xstar}=-\nabla f(\xstar)-\nabla_x g(\xstar,\lstar)=-\nabla f(\xstar)+J_h(x)^\top\lstar=0$.
\end{proof}

{As a consequence of the strong convexity of $f(x)$,} from \cite[Lemma 1]{qu19}, there exists a symmetric, positive definite $B=B(x)\in\R^{n,n}$ such that
\begin{equation}\label{lem:qu}
 \nabla f(x)-\nabla f(\xstar) = B(x-\xstar).
\end{equation}
%We write $B=B(x)$ for brevity.

\begin{theorem}[Global exponential convergence]\label{theo:pi}
%
%Let assumptions \ref{ass:sc} and \ref{ass:aff} hold.
%For any $\kp>0$, given a constant
%\begin{equation}\label{rho}\rho\geq \frac{\kp(\beta_2-\beta_1)^2}{2\beta_1},\end{equation}
Let us assume that $C$ is full rank and there exists $0<\underline{c}\leq \overline{c}$ such that
\begin{equation}\label{CCposdef}
 \underline{c} I \preceq CC^\top \preceq \overline{c} I.
\end{equation}
Let $\rho<\overline{c}^{-1}$. Then, there exist real positive constants $\alpha_1$ and $\alpha_2$ such that
\begin{equation}
 \|x(t)-\xstar\|_2\leq \alpha_1e^{-\frac{1}{2}\mu t},~~~\|\lambda(t)-\lstar\|_2\leq \alpha_2e^{-\frac{1}{2}\mu t}
 \end{equation}
where
\begin{equation}\label{tau}
\mu\leq \min\left\{ \frac{1}{2}\kp \underline{c},\frac{2\ki \underline{g}-\kp\overline{g}^2}{\ki} \right\}
\end{equation}
and $0<\underline{g}\leq \overline{g}$ are assessed in the proof.
\end{theorem}
\begin{proof}
We define the candidate Lyapunov function
\begin{equation}\label{def:lyap}
V\big(z(t)\big)=\big(z(t)-\zstar\big)^\top P \big(z(t)-\zstar\big)
\end{equation}
where
\begin{equation}\label{def:P}
P:=\left(\begin{array}{cc}
                    \sigma I_n&0\\
                    0&I_m
                   \end{array}\right)\in\R^{m+n,m+n}
            \end{equation}
for some $\sigma>0$. If
\begin{equation}\label{lyap}
 \dot{V}\big(z(t)\big)\leq -\mu V\big(z(t)\big)
\end{equation}
then the theorem statement holds. Therefore, let us study the conditions that guarantee \eqref{lyap}.

Let us consider the diagonal matrix $\Gamma = \Gamma(z)\in [0,1]^{m,m}$ as defined in Lemma 3 in \cite{qu19}.
Since $\nabla_x\Lag(\xstar,\lstar)=\nabla_{\lambda}\Lag(\xstar,\lstar)=0$, $J_h(x)=C$ and by using \eqref{lem:qu},
\begin{equation}\label{gradx}
 \begin{split}
  \dot{x}&=-\nabla_x\Lag(x,\lambda)=-\nabla_x\Lag(x,\lambda)+\nabla_x\Lag(\xstar,\lstar)\\
  &=-B(x-\xstar)- \rho C^{\top}\Gamma C(x-\xstar) -C^{\top}\Gamma(\lambda-\lstar)
 \end{split}
\end{equation}
as obtained for PDGD, see Sec. III-B in \cite{qu19} for details.
Furthermore,
\begin{equation}\label{gradl}
 \begin{split}
 \dot{\lambda}=&\nabla_{\lambda}\Lag(x,\lambda)=\nabla_{\lambda}\Lag(x,\lambda)-\nabla_{\lambda}\Lag(\xstar,\lstar)\\
  &=\ki \Gamma C (x-\xstar)+\frac{\ki}{\rho}(\Gamma - I)(\lambda-\lstar)+\kp C \dot{x}.
 \end{split}
\end{equation}
Equations \eqref{gradx}-\eqref{gradl} represent \eqref{PI} in a ``linear'' form.
Let $G_1:= B+\rho C^\top\Gamma C$ and $G_2:= C^\top \Gamma$. Since $B$ is positive definite, $G_1$ is positive definite; let {$\underline{g} I \preceq G_1^\top \preceq \overline{g} I.$}
 Then, we can rewrite \eqref{gradx}-\eqref{gradl} in a matrix form
\begin{equation}
\dot{z}(t)=G \big(z(t)-\zstar\big)
\end{equation}
where
\begin{equation}
G:=\left(\begin{array}{cc}
                    -G_1& -G_2\\
                    \ki G_2^\top-\kp C G_1&\frac{\ki}{\rho}(\Gamma-I)-\kp C G_2
                   \end{array}\right)
            \end{equation}
Since
\begin{equation}
\begin{split}
 \dot{V}\big(z(t)\big)&=\dot{z}(t)^{\top}P\big(z(t)-\zstar\big)+\big(z(t)-\zstar\big)^\top P\dot{z}(t)\\
 &=\big(z(t)-\zstar\big)^\top\big(G^\top P+P G\big)\big(z(t)-\zstar\big)
\end{split}
\end{equation}
a sufficient condition for  $\dot{V}\big(z(t)\big)\leq -\mu V\big(z(t)\big)$, see \eqref{lyap},
is
\begin{equation}\label{cond}
-G^{\top} P -P G-\mu P \succeq 0.
\end{equation}
As a consequence, our next goal is to provide sufficient conditions for \eqref{cond}.
%
% We compute
% %
% \begin{equation}\label{pg}
% PG=\left(\begin{array}{cc}
%                     -\rho B&-\rho C^\top\\
%                     K_i C -K_pCB&-K_p CC^\top
%                    \end{array}\right)
% \end{equation}
% while $G^{\top} P=(PG)^\top$. Hence,
%
We have
\begin{equation}\label{gg}
\begin{split}
 &-G^{\top} P -P G-\mu P= \left(\begin{array}{cc}
                   Q_1& Q_2\\
                    Q_2^\top&Q_3
                   \end{array}\right).
                   \end{split}
\end{equation}
where
\begin{equation}
 \begin{split}
  Q_1 &= 2\sigma G_1-\sigma \mu I_n\\
  Q_2  & =  (\sigma-\ki)G_2+\kp G_1^\top C^\top\\
  Q_3 & = \kp C G_2+\kp G_2^\top C^\top+2\frac{\ki}{\rho}(I-\Gamma)-\mu I_m
 \end{split}
\end{equation}

If $\ki\geq \kp$, by applying \cite[Lemma 6]{qu19} for $\frac{1}{\rho}>\overline{c}$,
\begin{equation}
 \kp C G_2+\kp G_2^\top C^\top+2\frac{\ki}{\rho}(I-\Gamma)\succeq \frac{3}{2}\kp C C^\top.
\end{equation}
Thus,
\begin{equation}\label{gg3}
Q_3 \succeq  \frac{3}{2}\kp C C^\top-\mu I_m \succeq \kp C C^\top
\end{equation}
where the last step is a consequence of the assumption $\mu\leq \frac{1}{2}\kp \underline{c} $.

To simplify the computations, we set  $\ki=\sigma$. Then,
\begin{equation}
   Q_2   =  \kp G_1^\top C^\top.
\end{equation}
In conclusion,
\begin{equation}\label{gge}
\begin{split}
 &-G^{\top} P -P G-\mu P\succeq \left(\begin{array}{cc}
                  2\sigma G_1-\sigma \mu I_n& \kp G_1^\top C^\top\\
                    \left[\kp G_1^\top C^\top\right]^\top& \kp C C^\top
                   \end{array}\right).
                   \end{split}
\end{equation}

% Let us define
% \begin{equation}\label{defQ}
%  Q:=(\rho-K_i)I +K_pB\in\R^{n,n}.
% \end{equation}
%
Since $CC^\top \succ 0$ is invertible from \eqref{CCposdef}, we can apply the Schur complement argument: the matrix in \eqref{gge}
is positive semidefinite if and only if
\begin{equation}\label{schur}
2\sigma G_1-\sigma \mu I_n- \kp G_1^\top C^\top \left(\kp C C^\top\right)^{-1} \kp C G_1 \succeq 0.
\end{equation}
Since $CC^\top$ is invertible, then $C^{\top}(CC^\top)^{-1}C\preceq I$. Therefore, a sufficient condition for \eqref{schur} is
\begin{equation}\label{sufficient}
2\sigma G_1-\sigma \mu I_n-  \kp G_1^\top G_1 \succeq 0.
\end{equation}
Finally,
\eqref{sufficient} holds if
\begin{equation}
 %2\sigma \underline{g}-\sigma \mu I_n-  \kp \overline{g}^2\geq 0
 2\sigma \underline{g} {-\sigma \mu} -  \kp \overline{g}^2\geq 0
\end{equation}
which is equivalent to
\begin{equation}
 \mu \leq \frac{2\ki \underline{g}-\kp\overline{g}^2}{\ki},
\end{equation}
which completes the proof.
\end{proof}
\begin{remark}
In the proof, we set {$\ki=\sigma$} to simplify the computations. Other choices may enhance the convergence rate.
\end{remark}
\begin{remark}
A theoretical comparison of the convergence rates $\mu$ and $\tau_{ineq}$ in \cite{qu19} is challenging because $\tau_{ineq}$ depends on many constants that are not easy to assess. Conversely, we can estimate $\mu$ straightforwardly, given some knowledge of the given optimization problem.
\end{remark}
\begin{remark}
The proof of Theorem \ref{theo:pi} {of this manuscript is more straightforward} than the proof of Theorem 2 in \cite{qu19} because the diagonal form of the Lyapunov function \eqref{lyap} reduces the computations {when compared to} the Lyapunov function with cross terms in \cite{qu19}.
\end{remark}

%{Aggiungere considerazioni su come scegliere $\ki,\kp$? Chiesto da 2 revisori su 3.}
{Theorem \ref{theo:pi} also suggests some insights on the selection of $\ki,\kp$. We notice that the ratio $\kp$ must be kept small to avoid reducing the convergence rate, while $\ki$ plays the same role of $\eta$ in \eqref{PDGD}.}

\subsection{Example on the convergence rate of PI and PDGD}
To conclude this section, we report an example to compare the convergence rates of PDGD in \eqref{PDGD} and PI in \eqref{PI}. %In particular, in this example PI is faster than PDGD. %In particular, we show that in some cases PI has a faster convergence rate than PDGD.

We consider the quadratic optimization problem
\begin{equation}
\begin{split}
&\min_{x\in\R} \frac{1}{2} w x^2\\
&\text{s.t.}\\
&x\leq 0
\end{split}
\end{equation}
where $w>0$.

Since the cost function is quadratic and the constraints are linear, the closed-loop dynamics with the PI control corresponds to a switched linear time-invariant system with the following two modes:

\begin{equation}\label{modes}
\begin{pmatrix}
                    \dot{x}\\
                    \dot{\lambda}\\
                   \end{pmatrix}
                   = \mathrm{A_1}
                  \begin{pmatrix}
                    x\\
                    \lambda\\
                   \end{pmatrix},~~~\begin{pmatrix}
                    \dot{x}\\
                    \dot{\lambda}\\
                   \end{pmatrix}
                   = \mathrm{A_2}
                  \begin{pmatrix}
                    x\\
                    \lambda\\
                   \end{pmatrix}              
\end{equation}
where
\begin{equation}\label{PI_modes}
\begin{split}
    \mathrm{A_1}&=\begin{pmatrix}
       -w-\rho&-1\\
    \ki-\kp(w+\rho)&-\kp 
    \end{pmatrix}\\
    \mathrm{A_2}&=\begin{pmatrix}
       -w&0\\
    -w\kp&-\frac{\ki}{\rho} 
    \end{pmatrix}.
    \end{split}
\end{equation}
In the first mode, the control through $\lambda$ is active; in the second mode, $x(t)$ satisfies the constraints and evolves according to the derivative of the cost function.

Similarly, {for PDGD, the dynamics correspond }to a switched linear time-invariant system of kind \eqref{modes} with two modes described by 
\begin{equation}\label{PDGD_modes}
    \mathrm{\tilde{A}_1}=\begin{pmatrix}
       -w-\rho&-1\\
    \eta&0 
    \end{pmatrix},~~~ \mathrm{\tilde{A}_2}=\begin{pmatrix}
       -w&0\\
    0&-\frac{\eta}{\rho} 
    \end{pmatrix}
\end{equation}
which correspond to $\mathrm{A_1}$ and $\mathrm{A_2}$ with $\kp=0$ and $\ki=\eta$.
Now, we notice that $\mathrm{A_2}$ and  $\mathrm{\tilde{A}_2}$ have the same eigenvalues for $\ki=\eta$, i.e., $-w$ and $-\frac{\ki}{\rho}$. Therefore, {PI and PDGD enjoy the same convergence rate in the second mode.}
Concerning the first mode, the eigenvalues of $\mathrm{A_1}$ are 
\begin{equation}
\frac{-\kp-w-\rho\pm\sqrt{(\kp+w+\rho)^2-4\ki}}{2}.
\end{equation}
Therefore, { if $\ki$ is }sufficiently large, the eigenvalues are complex conjugate, and the convergence rate is $\kp+w+\rho$.

On the other hand, the eigenvalues of $\mathrm{\tilde{A}_1}$ are 
\begin{equation}
\frac{-w-\rho\pm\sqrt{(w+\rho)^2-4\eta}}{2}
\end{equation}
with {the }best convergence rate equal to $w+\rho$. In conclusion, PI has a better convergence rate than PDGD, provided {a suitable $\kp$ is selected.} In particular, for PDGD, {increasing the convergence rate beyond $w+\rho$ is impossible.} 
%{Da revisore 3. Non basta aumentare $\rho$ in PDGD?}

{We remark that although it is always possible to increase $\rho$ both in \eqref{PI} and \eqref{PDGD}, this may cause numerical issues during the integration of the differential equations.}
\section{Numerical results}\label{sec:NR}
In this section, we illustrate two numerical simulations to support the effectiveness of the proposed algorithm \eqref{PI} compared to PDGD.

\subsection{Example 1: Quadratic programming}
In this simulation, we consider a strongly convex quadratic programming (QP) in the following form:
\begin{equation}\begin{aligned}
    x^{\star} &= \arg\min_{x\in\R^n} \frac{1}{2} x^\top \left( I + W^\top W \right) x + b^\top x \\
    &\text{s.t.} \\
    &C x -d \leq 0
\end{aligned}\end{equation}
where $W \in \R ^{n,n},b \in \R ^n, C \in \R ^{m,n},d \in \R ^m$ are randomly generated vectors or matrices, with independent, normally distributed components. We set $n=50$ and $m=45$.

We solve the optimization problem using the proposed algorithm \eqref{PI} and PDGD. We integrate the differential equations \eqref{PDGD} and \eqref{PI} in the time interval $[0,30]\,\rm{s}$ through the \textit{ode45} MATLAB command to select the discretization step size optimally. We set $\ki=\eta=1$ and $\kp=-0.7$.

\begin{figure}
    \centering
    \includegraphics[width=\linewidth]{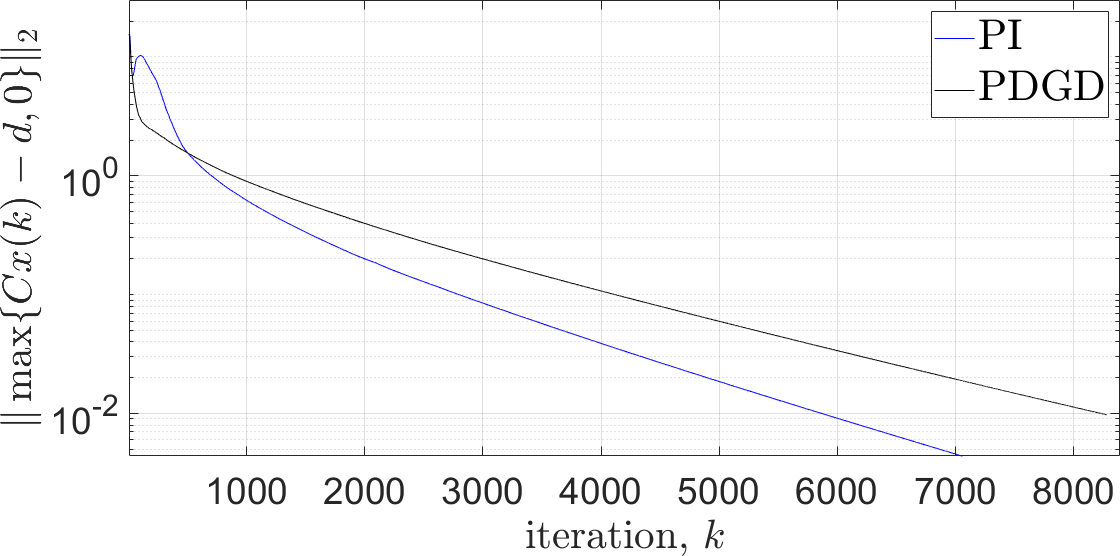}
    \caption{Example 1: constraints violation.}
    \label{fig:qp_cns}
\end{figure}

Fig. \ref{fig:qp_cns} shows the time evolution of $\|\max{Cx(k)-d,0}\|_2$, where $k$ is the current iteration. This metric represents the violation of the constraints and is equal to zero when the state $x$ satisfies the constraints.

\begin{figure}
    \centering
    \includegraphics[width=\linewidth]{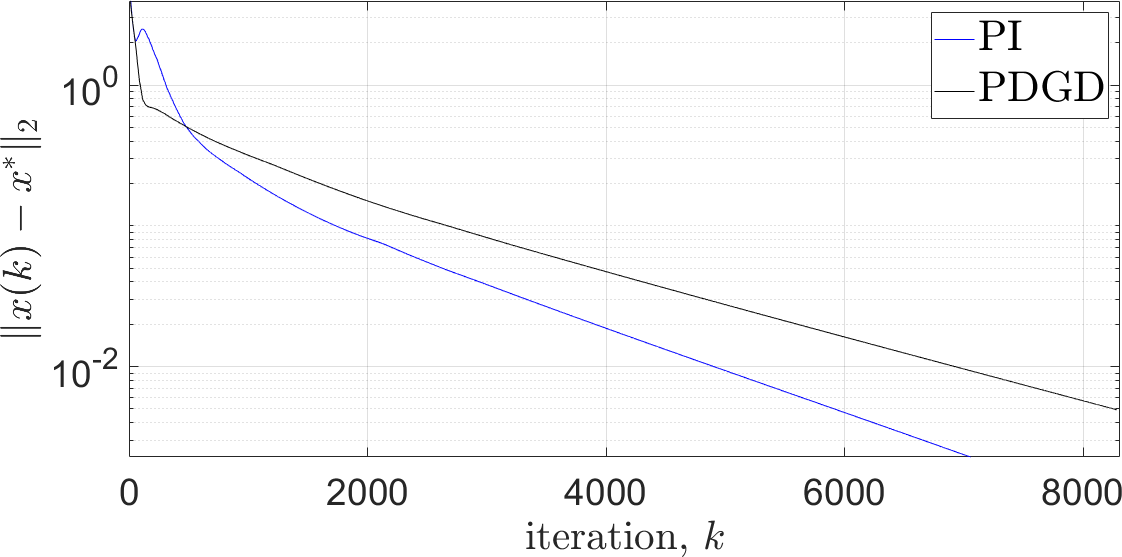}
    \caption{Example 1: distance from the optimum.}
    \label{fig:qp_opt}
\end{figure}

Fig. \ref{fig:qp_opt} shows the $\ell_2$ distance from the global optimum $x^{\star}$, computed through the MATLAB package CVX \cite{cvx}.


We perform 100 random runs with different realizations of $W, b, C, d$. PI requires fewer iterations than PDGD in all the runs and, {consequently}, less computational time. In Table \ref{tab:es1}, we report some statistics that show the enhanced convergence speed of PI with respect to PDGD.
%
%In particular, PI requires $6903.2 \pm 312.07$ iterations corresponding to a computation time of $(1.02\cdot10^{-1} \pm 4.74\cdot10^{-3}) \,\rm{s}$, while PDGD requires $8128.3 \pm 491.94$ iterations corresponding to a computation time of $(1.24\cdot10^{-1} \pm 8.04\cdot10^{-3}) \,\rm{s}$.

Fig. \ref{fig:qp_cns} and Fig. \ref{fig:qp_opt}, obtained from one randomly selected run, show us that the proposed approach converges more quickly than PDGD, either in terms of fulfilment of the constraints and achievement of the minimum. 

\begin{table}[h]
    \centering
    {\small{}
    \begin{tabular}{l|c c c}
         & mean & standard deviation & worst case \\
         \hline
        $N$ PDGD & $8128.3$ & $491.9$ & $ 9361$ \\
        $N$ PI & $6903.2 $ & $312.1$ & $7613$ \\
        $T$ PDGD & $\SI{1.23e-01}{}$  & $\SI{8.04e-03}{}$ & $\SI{1.53e-01}{}$\\
        $T$ PI & $\SI{1.02e-01}{}$ & $\SI{4.74e-03}{}$ & $\SI{1.15e-01}{}$ \\
        %\hline 
        %$E$ PDGD & $\SI{5.86e-03}{}$ & $\SI{1.03e-02}{}$ & $\SI{4.61e-02}{}$ \\
        %$E$ PI & $\SI{3.83e-03}{}$ & $\SI{7.66e-03}{}$ & $\SI{3.85e-02}{}$ \\
    \end{tabular}}
    \caption{Example 1: results over 100 random runs. $N$ is the required numbers of iteration; $T$ the computational time (in seconds).}
    \label{tab:es1}
\end{table}

\subsection{Example 2: Linear system identification}
We apply our approach to a problem of system identification. We consider the problem of identifying an unknown stable linear dynamic system $H(z)$ using uncertain input-output measurements $\{u_k,\Tilde{y}_k\}$, for $k\in\{1,\dots, N\}$, where $\Tilde{y}_k = y_k+\eta_k$, $y_k$ is the $k$-th noise-free output sample and $\eta_k$ is the $k$-th sample of the noise sequence.% We assume that all samples of the noise $\eta_k$ are unknown but bounded in magnitude by a given constant $\Delta_\eta$, i.e., $\Vert y_k-\Tilde{y}_k \Vert \leq \Delta_\eta$ for all $k\in\{1,\dots, N\}$.

To perform the identification, we select a linear-in-the-parameters model structure $\Tilde{H}(z,\theta)$ of the form 
\begin{equation}
    \Tilde{H}(z) = \sum_{i=1}^P \theta_i B_i(z),
\end{equation}
which is a standard choice in the context of system identification, and commonly considered choices of $B_i(z)$ are Laguerre or Kautz filters; see, e.g., \cite{van1995system},\cite{wahlberg1996approximation}. For the sake of simplicity, in this example, we select $B_i(z)$ to be first-order transfer functions with poles uniformly distributed in $[-0.9, 0.9]$.
More precisely, we select 
\begin{equation}
    B_i(z) = \frac{z}{z-p_i}, \quad p_i \in \{-0.9,-0.85,\dots,0.9\}. %\bigcup \bigg\{ \frac{z}{(z-p_i)(z-p_i^*)}, \\p_i=0.1 e^{i\gamma},\dots,0.9 e^{i\gamma}, \gamma = \pi/4,\dots,7\pi/4\bigg\}.
\end{equation}

To generate time-domain data, we simulate the randomly selected system
\begin{equation}
H(z) = \frac{-0.4 z^2 + 0.32 z + 0.26}{z^3 - 1.9 z^2 + 1.21 z - 0.259}
\end{equation}
excited by a random input uniformly distributed in $[0,1]$. We corrupt the output data with normally distributed noise with variance $\sigma^2_\eta = 0.1$.

We look for the value of the parameter $\theta$ that minimizes the $\ell_\infty$-norm of the simulation error
\begin{equation}\label{eq:id_probl}\begin{aligned}
    &\theta^* = \arg\min_{\theta \in \R^P}  \norm{y_k(\theta) - \Tilde{y}_k }_\infty.
\end{aligned} \end{equation}
By adding a slack variable $\Delta \in \R$, we recast problem \eqref{eq:id_probl} to the following linear programming problem: 
\begin{equation}\label{eq:id_probl2}\begin{aligned}
    &\theta^*, \Delta^* = \arg\min_{\theta \in \R^P, \Delta \in \R}  \Delta \\
    &\qquad \text{s.t.} \\
    &\quad -\Delta \leq Z_k \theta - \Tilde{y}_k \leq \Delta, \quad k \in \{1,\dots,N\}
\end{aligned} \end{equation}
where $Z_k = [z_1(k),\dots,z_P(k)] \in \R^{P}$, $z_i(k) = B_i(q^{-1})u(k)$.

We solve the optimization problem using the proposed PI algorithm in \eqref{PI} and PDGD in \eqref{PDGD}. We integrate the differential equations \eqref{PDGD} and \eqref{PI} in the time interval $[0,1000]\,\rm{s}$ through the \textit{ode23} MATLAB command to select the discretization step size optimally. We set $\ki=\nu=1$ and $\kp=-0.5$. 

\begin{figure}
    \centering
    \includegraphics[width=\linewidth]{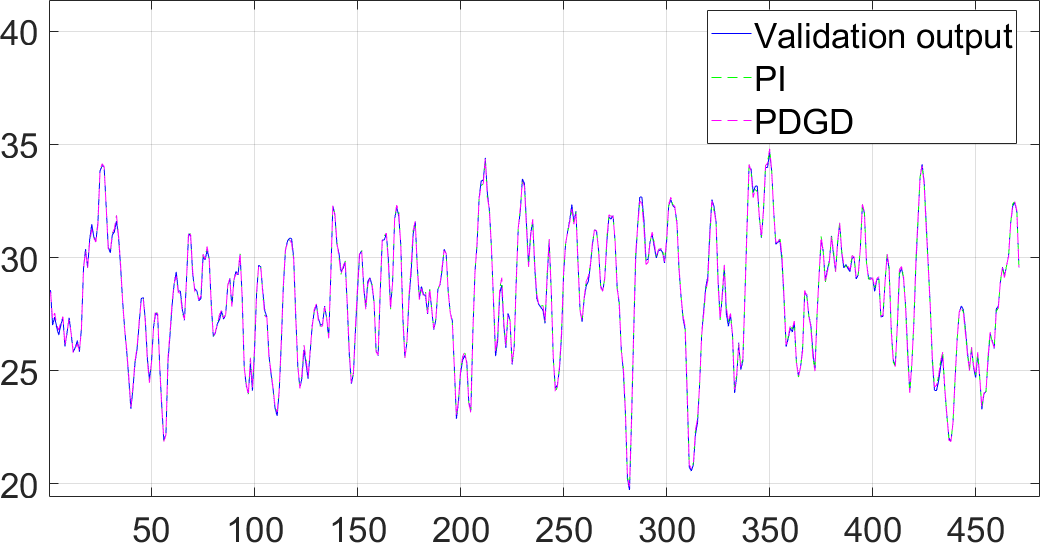}
    \caption{Example 2: validation of identified models}
    \label{fig:es02}
\end{figure}

In Fig. \ref{fig:es02}, we compare the outputs of the true model and the ones estimated using the PDGD and PI algorithms on data not used for identification. The outputs of the three models are almost exactly overlapped. We also evaluate the validation performances of the two algorithms in terms of the FIT index, defined as
\begin{equation}
    FIT = 100\left(1-\sqrt{\frac{\norm{y^{val} - \hat{y}^{val}}^2_2}{\norm{y^{val} - m_y^{val}}^2_2}}\right)
\end{equation}
where $y^{val}$ is the true output, $m_y^{val}=\frac{1}{N}\sum_{k=1}^Ny^{val}_k$ and $\hat{y}^{val}$ is the response of the identified model. We obtain the same value of $\mbox{FIT}=98.5\%$ for both the identified models. Fig. \ref{fig:es02} and the computed FIT values show that both algorithms converge to the optimal solution as expected from the theory. We report the comparison of the two algorithms in terms of computational effort in Table \ref{tab:es2}, where we show the required number of iterations and computational time. Such results show that the PI algorithm is about two times faster than PDGD.

%{Beyond the proposed method and PDGD, we solve this problem also through the CVX MATLAB toolbox \cite{cvx} that by default relies on the SDPT3 solver \cite{tut03}. We notice that SDPT3 fails to converge to the optimal solution due to the bad conditioning of the problem, i.e., $\text{cond}(Z Z^\top) \approx 10^{20}$. For the same reason, according to $\tau_{\rm ineq}$ in \cite{qu19} and $\mu$ in \eqref{tau}, the convergence rate of both PDGD and PI is slow. In Table \ref{tab:es2}, we compare PI and PDGD in terms of the required number of iterations and computational times.}

%As to the convergence rate, we report the number of iterations and computational times in Table \ref{tab:es2}.

%PDGD requires about $1\cdot10^{6}$ iterations and $\SI{679.46}{\second}$ to achieve convergence, while PI requires about $514\cdot10^3$ iterations and $\SI{233.09}{\second}$.

\begin{table}[]
    \centering
    {\small{}
    \begin{tabular}{l|c c}
         & $N$ & $T$  \\
         \hline
         PDGD & $10^6$ & $679.46$  \\
         PI & $5\times 10^5 $ & $233.09$\\
        %\hline 
        %$E$ PDGD & $\SI{5.86e-03}{}$ & $\SI{1.03e-02}{}$ & $\SI{4.61e-02}{}$ \\
        %$E$ PI & $\SI{3.83e-03}{}$ & $\SI{7.66e-03}{}$ & $\SI{3.85e-02}{}$ \\
    \end{tabular}}
    \caption{Example 2: $N$ is the required numbers of iteration; $T$ the computational time (in seconds).}
    \label{tab:es2}
\end{table}

%Beyond the proposed method and PDGD, we solve this problem also through the CVX MATLAB toolbox \cite{cvx} that relies on the SDPT3 solver \cite{tut03}. Interestingly, SDPT3 fails to converge to the optimal solution due to the bad conditioning of the problem (i.e., $\text{cond}(Z^\top Z) \approx 10^{19}$. On the other hand, both PDGD and PI are implemented through a Simulink scheme to exploit the automatic integration step-size selection, and we select $\ki=1$ for both and $\kp=0.1$ for PI. Both converge quickly to the optimal solution. In particular, PI requires $300$ steps and $0.9488$ seconds to converge, while PDGD requires $326$ steps corresponding to $1.177$ seconds.

\section{Conclusion}\label{sec:CON}

This paper proposes a novel continuous-time algorithm to solve smooth, strongly convex optimization problems with inequality constraints. As for the primal-dual gradient dynamics, the proposed algorithm consists of a dynamic system in the optimization variable and the Lagrange multipliers of the problem. By elaborating on the feedback PI control approach proposed for the equality-constrained case in \cite{cmo24}, we develop a variant of primal-dual gradient dynamics in which an additional term adjusts the dynamics of the Lagrange multipliers and enhances the convergence speed. We prove that the proposed method is globally exponentially convergent. Finally, examples and numerical simulations show its effectiveness and velocity concerning the PDGD algorithm. 
{ Current work envisages the relaxation of affinity and smoothness requirements in the considered model.}
%{ Future works will envisage relaxing the requirement that $g(x)$ are affine functions and that $f(x)$ is smooth and strongly convex.}
%Finally, we notice that it is possible to combine the proposed approach with the PI control method in \cite{cmo24} to deal with equality and inequality constraints simultaneously. 

\bibliographystyle{IEEEtran}
\bibliography{sophie}

\end{document}